\documentclass{ifacconf}

\usepackage{graphicx}      
\usepackage{natbib}        
\newcommand{\He}{\text{He}}
\newcommand{\wh}{\hat{w}}
\newcommand{\zh}{\hat{z}}
\newcommand{\uh}{\hat{u}}
\newcommand{\vh}{\hat{v}}
\newcommand{\nh}{n^s}

\usepackage{amssymb}									
\usepackage{amsmath}
\usepackage{amsfonts}
\usepackage{mathrsfs} 								
\usepackage[tikz]{bclogo}
\usepackage{color,colortbl}		
\usepackage{arydshln}

\definecolor{xgray}{rgb}{0.75, 0.75, 0.75}
\newcommand{\hdl}{\hdashline[2pt/1pt]}
\newcommand{\hgl}{\noalign{\global\arrayrulewidth=0.6pt}\arrayrulecolor{xgray}\hline\noalign{\global\arrayrulewidth=0.4pt}\arrayrulecolor{black}}
\newcommand{\vgl}{\color{xgray}\vrule width0.6pt}

\let\originalleft\left
\let\originalright\right
\renewcommand{\left}{\mathopen{}\mathclose\bgroup\originalleft}
\renewcommand{\right}{\aftergroup\egroup\originalright}

\newcommand{\R}{\mathbb{R}}
\newcommand{\Sb}{\mathbb{S}}
\newcommand{\Ahb}{\mathbb{A}}
\newcommand{\Bhb}{\mathbb{B}}
\newcommand{\Chb}{\mathbb{C}}

\newcommand{\Ac}{\mathcal A}
\newcommand{\Bc}{\mathcal B}
\newcommand{\Cc}{\mathcal C}
\newcommand{\Dc}{\mathcal D}
\newcommand{\Kc}{K}
\newcommand{\Lc}{\mathcal L}
\newcommand{\Pc}{\mathcal P}
\newcommand{\Xc}{\mathcal X}
\newcommand{\Yc}{\mathcal Y}
\newcommand{\Zc}{\mathcal Z}
\newcommand{\AB}{\boldsymbol{A}}
\newcommand{\BB}{\boldsymbol{B}}
\newcommand{\CB}{\boldsymbol{C}}
\newcommand{\DB}{\boldsymbol{D}}
\newcommand{\NB}{\boldsymbol{N}}
\newcommand{\XB}{\boldsymbol{X}}
\newcommand{\KB}{\boldsymbol{K}}
\newcommand{\LB}{\boldsymbol{L}}
\newcommand{\MB}{\boldsymbol{M}}

\newcommand{\Qa}{\tilde{Q}}
\newcommand{\Sa}{\tilde{S}}
\newcommand{\Ra}{\tilde{R}}
\newcommand{\Ah}{\hat{A}}
\newcommand{\Bh}{\hat{B}}
\newcommand{\Ch}{\hat{C}}

\newcommand{\Ms}{\bar{M}}
\newcommand{\As}{\bar{A}}
\newcommand{\Bs}{\bar{B}}
\newcommand{\Css}{\bar{C}}

\newcommand{\Ks}{\bar{K}}
\newcommand{\Ls}{\bar{L}}
\newcommand{\al}{\alpha}
\newcommand{\ga}{\gamma}
\newcommand{\de}{\delta}
\newcommand{\De}{\Delta}

\newcommand{\mat}[2]{\left(\begin{array}{#1}#2\end{array}\right)}
\newcommand{\sy}{(\ast)^T}
\newcommand{\diag}{\operatorname{diag}}
\newcommand{\tr}{\operatorname{tr}}
\newcommand{\col}{\operatorname{col}}

\newcommand{\smat}[1]{\left(\begin{smallmatrix}#1\end{smallmatrix}\right)}
\newcommand{\cl}{\prec}
\newcommand{\cg}{\succ}
\renewcommand{\t}{\tilde}
\newcommand{\ti}{\times}
\newcommand{\Hz}{\mathcal{H}_2}

\usetikzlibrary{matrix}

\usepgflibrary{fpu}
\usetikzlibrary{positioning,calc,arrows,shapes,shadows,matrix,backgrounds}

\tikzset{
	auto,
	poi/.style={
		minimum size=0,
		inner sep=0
	},
	sys/.style 2 args={
		rectangle,
		draw,
		rounded corners,
		drop shadow,
		draw=black,
		top color=black!20,bottom color=black!0,
		minimum height=#2,
		minimum width=#1,
		inner sep=\dn},
	syse/.style 2 args={
		rectangle,
		draw=none,
		rounded corners,
		minimum height=#2,
		minimum width=#1,
		inner sep=\dn},
	nod/.style={
		circle,
		draw,
		fill=white,
		minimum size=5ex
	},
	sum/.style={circle,draw,draw=black,inner sep=0mm,minimum size=2mm,fill=white,
		draw=black!100,top color=black!20,bottom color=black!0},
	jun/.style={circle,draw,draw=black,inner sep=0mm,minimum size=0mm},
	>={latex},
	every path/.style={rounded corners},
	lin/.style={color=black,draw,->},
	ju/.style={inner sep=0mm,minimum size=5*\dn},
	sysr/.style 2 args={
		circle,
		draw,
		drop shadow,
		draw=black,
		top color=black!20,bottom color=black!0,
		minimum height=#2,
		minimum width=#1,
		inner sep=\dn},
}

\def\dn{1ex}

\tikzstyle{syr4}=[sysr={3*\dn}{3*\dn}] 
\tikzstyle{syr6}=[sysr={6*\dn}{6*\dn}] 
\tikzstyle{syr7}=[sysr={7*\dn}{7*\dn}] 

\usepackage[absolute]{textpos} 

\begin{document}
\begin{frontmatter}

\title{Lifting to Passivity for $\Hz$-Gain-Scheduling Synthesis with Full Block Scalings\thanksref{footnoteinfo}}

\thanks[footnoteinfo]{
Funded by Deutsche Forschungsgemeinschaft (DFG, German Research Foundation) under Germany's Excellence Strategy - EXC 2075 - 390740016.}

\author[First]{Christian A. R\"osinger and Carsten W. Scherer}

\address[First]{Department of Mathematics, University of Stuttgart, Germany,\\\hspace{-3.5ex} e-mail: \{christian.roesinger,carsten.scherer\}@mathematik.uni-stuttgart.de}

\begin{abstract}
We focus on the $\Hz$-gain-scheduling synthesis problem for time-varying parametric scheduling blocks with scalings. Recently, we have presented a solution of this problem for $D$- and positive real scalings
by guaranteeing finiteness of the $\Hz$-norm for the closed-loop system with suitable linear fractional plant and controller representations. In order to reduce conservatism, we extend these methods to full block scalings by designing a triangular scheduling function and by introducing a new lifting technique for gain-scheduled synthesis that enables convexification.
\end{abstract}


\end{frontmatter}

\begin{textblock}{12}(1.75, 15.75)
\fbox{
\begin{minipage}{\textwidth}
\small\textcopyright~ 2020 IFAC. This work has been published in IFAC-PapersOnline under a Creative Commons Licence CC-BY-NC-ND. DOI: 10.1016/j.ifacol.2020.12.570
\end{minipage}
}
\end{textblock}

\section{Introduction}

The design of linear parametrically-varying (LPV) systems is widely spread over the control literature and can be roughly divided into two classes. On the one hand, parameter-dependent Lyapunov functions, as in \cite{becker1995}, \cite{wu1996}, \cite{apkarian1997}, \cite{wu2005}, \cite{souza2006}, and \cite{sato2011}, are used for synthesis with linear matrix inequalities (LMIs) by approximating the parameter space of the scheduling variable. On the other hand, the so-called scaling approach can directly handle rational parameter dependence, as in \cite{packard1994}, \cite{apkarian1995} for $D$-scalings, \cite{helmersson98} for positive real-scalings, \cite{scorletti1998} for $D/G$-scalings, and \cite{scherer2000}, \cite{veenman2014} for the least conservative full block scalings. These approaches are as well of interest because of their link to distributed controller design (see \cite{langbort2004}) and their flexibility for handling more complex scheduling blocks such as delays as considered in \cite{roesinger2019b}.

\begin{figure}
	\vspace{-1ex}
	\centering
	\includegraphics[width=0.2\textwidth]{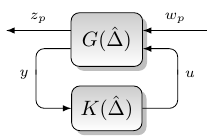}\label{GS2}
	\caption{Feedback-loop for gain-scheduling.}
\end{figure}

In this work, we look at the concrete configuration in Fig.~\ref{GS2} which has shown to be well-suited for analysis and synthesis of LPV controllers (see \cite{packard1994}, \cite{apkarian1995}). For an uncertain plant $G(\hat{\De})$ with $\hat{\De}$ being an arbitrary fast time-varying matrix-valued parametric uncertainty, we employ constant full block scalings to synthesize a controller $K(\hat{\De})$ which achieves an $\Hz$-cost criterion imposed on $w_p\to z_p$. Concrete applications of LPV design with $\Hz$-performance guarantees are, e.g., the control of autonomous cars and helicopters in \cite{mustaki2019} and \cite{guerreiro2007}, respectively.\\
Recently, \cite{roesinger2019b} present the first scaling solution to this problem with $D$-scalings in case that the
uncertainty takes values in the unit disk or with positive-real scalings in case that the uncertainty is passive. Technically, this approach uses a convexifying transformation for controller and scaling parameters based on \cite{masubuchi1998}, \cite{scherer1997}, while suitable structured plant and controller descriptions guarantee well-posedness for the closed-loop $\Hz$-norm by design.
However, these results heavily rely on the particular structure of $D$- and positive real scalings and cannot be easily extended to the less conservative full block scalings.

As the main contribution of this work, we present a complete solution for the $\Hz$-gain scheduling problem with full block scalings in terms of LMIs. For this purpose, we introduce a new design approach based on what we call lifting to passivity.
This amounts to a loss-less embedding of the original synthesis problem into a passivity framework involving a suitable structural extension (or lifting) of the plant and the controller, and is the enabling factor for being able to convexify the problem through a transformation that operates on both the controller and the scaling parameters.
The use of a related passivation step has been beneficial already for a completely different objective in robustness analysis and synthesis involving integral quadratic constraints in \cite{veenman2013}, \cite{veenman2014}.
As a novel feature of this paper, we develop a systematic approach for using such a procedure in the context of gain-scheduled synthesis. 
As a further contribution, we reveal how suitable structured plant and controller representations can be exploited 
in our designs to render the $\Hz$-norm finite.

\textbf{Outline.}
After introducing the notation used in this work, Section~\ref{sec1} formulates the $\Hz$-gain scheduling problem under investigation, while Section~\ref{sec:lifting} presents the lifting design technique. The resulting specifically structured design problem is solved in Section~\ref{sec2}. Finally, a short example clarifies that our results are less conservative than those in \cite{roesinger2019b}.

\textbf{Notation.}
Let $\Sb^n$ denote the set of real symmetric matrices of dimension $n\ti n$.
For some matrices $M\in\R^{r\ti s}$ and $P\in\R^{r\ti r}$ we abbreviate $M^TPM$ by $\sy PM$ and $P+P^T$ by $\He(P)$ and denote by $\tr(P)$ the trace of $P$.
Matrix entries that can be inferred by symmetry are indicated by $\ast$.
We drop superscripts specifying partitions and dimensions of matrices if they are clear from the context.
Further, $I$ and $I_m$ denote identity matrices (with $m$ specifying the dimension if not clear from the context) and $\col(u_1,u_2):=\mat{cc}{u_1^T& u_2^T}^T$ is used for vectors. If $X,\,R,\,S$ and $A_{ij},\,B_i,\,C_j,\,D$ are some suitable matrices for $i,j=1,2$, we abbreviate
$$
\Lc\left(X,R,S,\smat{A_{ij}&B_i\\C_j&D}\right)
\,{:=}\,\sy\hspace{-0.3ex}\mat{c;{2pt/1pt}c;{2pt/1pt}c}{X&0&0\\\hdl 0&R&0\\\hdl 0&0&S}\hspace{-0.3ex}
\mat{ccc}{I&0&0\\ A_{11}&A_{12}&B_1\\\hdl 0&I&0\\ A_{21}&A_{22}&B_2 \\\hdl 0&0&I\\C_1&C_2&D}
$$
and refer to its left upper sub-block as
$$
\Lc_{\text{sub}}\left(X,R,S,\smat{A_{ij}\\C_j}\right):= \sy \mat{c;{2pt/1pt}c;{2pt/1pt}c}{X&0&0\\\hdl 0&R&0\\\hdl 0&0& S}
\mat{cc}{I&0\\A_{11}&A_{12}\\\hdl 0&I\\A_{21}&A_{22}\\\hdl 0&0\\C_1&C_2}.
$$

\section{Problem formulation}\label{sec1}

In the sequel, we introduce the $\Hz$-gain scheduling problem for full block scalings.

\subsection{Structured plant and controller representations}\label{sec21}

For some full block time-varying uncertainty $\hat{\De}$ taking values in some polytope, let
us consider the standard LPV configuration in Fig.~\ref{GS2} with a $\hat{\De}$-dependent LPV system $G(\hat{\De})$ and a corresponding controller $K(\hat{\De})$. 
To systemati\-cally guarantee finiteness of the closed-loop $\Hz$-norm, we use specifically structured linear fractional representations (LFRs) for $G(\hat{\De}),\,K(\hat{\De})$. Let $G(\hat{\De})$ be structured as in
\begin{equation}\label{decp1}
\mat{c}{\dot{x}\\\hline z_p\\ y}=\mat{c|cc}{A(\hat{\De})&B^p(\hat{\De})&B^u(\hat{\De})\\\hline
	C^p(\hat{\De}) &0 &D^u(\hat{\De})\\ C^y(\hat{\De}) &  D^y(\hat{\De})&D(\hat{\De})}\mat{c}{x\\\hline w_p\\ u},
\end{equation}
with $D(0)=0$, performance channel $w_p\to z_p$, control channel $u\to y$, and let us describe the controller $K(\hat{\De})$ by
\begin{equation}\label{controller}
\mat{c}{\dot{x}_c\\u}=\mat{cc}{A^c(\hat{\De})&B^c(\hat{\De})\\C^c(\hat{\De})&0}\mat{c}{x_c\\y}
\end{equation}
such that all $\hat{\De}$-dependent operator blocks in (\ref{decp1}), (\ref{controller}) are LFRs in $\hat{\De}$. Analogous to the approach for one repeated block in \cite{roesinger2019b}, the zero block structures in (\ref{decp1}), (\ref{controller}) guarantee that the performance channel in Fig.~\ref{GS2} has an identically vanishing direct feedthrough term.
Since $w_p\to z_p$ is zero in (\ref{decp1}), standard techniques for linear fractional transformations (LFTs) show that $G(\hat{\De})$ can be expressed as the LFR
\begin{eqnarray}\label{s00}
\begin{aligned}
&\mat{c}{\dot{x}\\\hline \hat{z}_1\\\hgl\hat{z}_2\\\hdl z_p\\ y}=\mat{c|c;{2pt/1pt}cc}{A_{11}&\Ah_{12}&B_1^p&B_1\\\hline \Ah_{21}&\Ah_{22}&\Bh_2^p&\Bh_2\\\hdl
	C_1^p &\Ch_2^p&D^p &D_1\\ C_1 & \Ch_2& D_2&D_3}\mat{c}{x\\\hline \hat{w}_1\\\hgl \hat{w}_2\\\hdl w_p\\ u}
=\\&\hspace{-1ex}{=}\,\renewcommand\arraystretch{1.2}
\mat{c|c!{\vgl}c;{2pt/1pt}cc}{A_{11}&\As_{12}&\As_{13}&B_1^p&B_1\\\hline\As_{21}&\As_{22}&0&0&\Bs_2\\\hgl\As_{31}&\As_{32}&\As_{33}&\Bs^p_3&\Bs_3\\\hdl C_1^p&\Css_2^p&0&0&D_1\\C_1&\Css_2&\Css_3&D_2&0}\hspace{-0.5ex}
\renewcommand\arraystretch{0.9}\mat{c}{x\\\hline \hat{w}_1\\\hgl \hat{w}_2\\\hdl w_p\\ u}
,\quad
\hat{w}=\De\hat{z}
\end{aligned}
\end{eqnarray}
with matrices $A_{11}\in\R^{n^s\ti n^s},\, B_1\in\R^{n^s\ti m},\, C_1\in\R^{k\ti n^s}$, as well as with
a structured uncertainty channel $\wh\to \zh$ for $\hat{z}:=\col(\hat{z}_1,\hat{z}_2)$ and $\hat{w}:=\col(\hat{w}_1,\hat{w}_2)$; the matrices associated to $\wh\to \zh$ are indicated with the symbols $\wedge$ or $-$ in (\ref{s00}).
W.l.o.g., the LFT manipulations can be always performed such that $\De=\diag(\hat{\De},\hat{\De})$ has a diagonal structure which is compatible with the partition of $\hat{A}_{22}$.\\
Since we only work with $\De$ in the sequel, we write $G(\De)$ for (\ref{s00}) and assume that $\De\in\mathbf{\De}$ where $\mathbf{\De}:=C([0,\infty),\mathbf{V})$ is the corresponding class of full block time-varying uncertainties for some given value set $\mathbf{V}=\text{Co}\{\De_1,\ldots,\De_N\}\ni0$ represented as the convex hull of finitely many real matrices $\De_i\in\R^{\uh\ti\vh}$. We hence consider (\ref{s00}) with $\De\in\mathbf{\De}$ as the precise mathematical description for (\ref{decp1}).\\
As the zero block structure for $K(\hat{\De})$ in (\ref{controller}) resembles that in (\ref{decp1}), the above LFT manipulations motivate to look at the following structured controller LFR
\begin{equation}\label{s12}
\begin{aligned}
&\mat{c}{\dot{x}_c\\\hline \hat{z}_{c,1}\\\hgl \hat{z}_{c,2}\\\hdl u}=
\mat{c|c;{2pt/1pt}c}{A_{11}^c&A_{12}^c&B^c_1\\\hline A_{21}^c&A_{22}^c&B^c_2\\ \hdl C^c_1&C^c_2&D^c}
\mat{c}{x_c\\\hline \hat{w}_{c,1}\\\hgl \hat{w}_{c,2}\\\hdl y}=\\
&\hspace{-0.5ex}{=}\,\renewcommand\arraystretch{1.2}\mat{c|c!{\vgl}c;{2pt/1pt}c}{A_{11}^c&\As_{12}^c&\As_{13}^c&B^c_1\\\hline \As_{21}^c&\As_{22}^c&0&0\\\hgl \As_{31}^c&\As_{32}^c&\As_{33}^c&\Bs_3^c\\\hdl C^c_1&\Css_2^c&0&0}
\renewcommand\arraystretch{1.0}\mat{c}{x_c\\\hline \hat{w}_{c,1}\\\hgl \hat{w}_{c,2}\\\hdl y},
\quad \hat{w}_c=\De_c(\De)\hat{z}_c
\end{aligned}
\end{equation}
with $\hat{z}_c:=\col(\hat{z}_{c,1}, \hat{z}_{c,2}),\, \hat{w}_c:=\col(\hat{w}_{c,1}, \hat{w}_{c,2})$ and the matrices $A_{11}^c\,{\in}\,\R^{n^c\ti n^c},\, B^c_1\,{\in}\,\R^{n^c\ti k},\, C^c_1\,{\in}\,\R^{m\ti n^c}$. We refer to (\ref{s12}) as $K(\De)$ in order to display the dependence on $\De$.
In order to have large enough flexibility in synthesis, we search for a lower block-triangular scheduling function
\begin{equation}\label{dec}
\De_c:\mathbf{V}\to\R^{r^c\ti r^c}\ \ \text{with}\ \ \De_c(V)\,{:=}\,\smat{\De^c_{11}(V)&0\\\De^c_{21}(V)&\De^c_{22}(V)}
\end{equation}
of partition $r^c:=r^c_1+r^c_2$. Indeed, for such a triangular $\De_c(.)$, the controller LFR (\ref{s12}) still ensures the structure in (\ref{controller}). Note that $\De_c(\De)$ might depend in a nonlinear fashion on $\De\in\mathbf{\De}$, while the choice of $r^c,\,n^c$ is part of the design problem.
The closed-loop system for the plant (\ref{s00}) interconnected with (\ref{s12})
is then given by
\begin{equation}\label{cl1}\renewcommand{\arraystretch}{1.0}\arraycolsep=1.0pt
\mat{c}{\dot{x}_e\\\hdl \hat{z}\\\hat{z}_c\\\hdl z_p}=\mat{c;{2pt/1pt}c;{2pt/1pt}c}{\hat{\Ac}_{11}&\hat{\Ac}_{12}&\hat{\Bc}_1\\\hdl \hat{\Ac}_{21}&\hat{\Ac}_{22}&\hat{\Bc}_2\\\hdl \hat{\Cc}_1&\hat{\Cc}_2&\hat{\Dc}}\mat{c}{x_e\\\hdl \hat{w}\\\hat{w}_c\\\hdl w_p},\, \mat{c}{\hat{w}\\\hat{w}_c}=\De_{ex}(\De)\mat{c}{\hat{z}\\\hat{z}_c}
\end{equation}
with extended state $x_e:=\col(x,x_c)$, extended scheduling block $\De_{ex}(V):=\smat{V&0\\0&\De_c(V)}$, and suitable closed-loop matrices $\hat{\Ac}_{ij},\,\hat{\Bc}_i,\,\hat{\Cc}_j,\,\hat{\Dc}$ for $i,j=1,2$.
\begin{defn}
	The controlled system (\ref{cl1}) is \textit{well-posed} if $I-\De_{ex}(V)\hat{\Ac}_{22}$ is non-singular for all $V\in\mathbf{V}$. It is \textit{stable} if there exist constants $K$ and $\al>0$ such that every solution of (\ref{cl1}) which is obtained for $w_p=0$ and any $\De\in\mathbf{\De}$ fulfills
	$$
	\|x_e(t)\|\leq Ke^{-\al(t-t_0)}\|x_e(0)\|\quad \text{for all}\quad t\geq 0.
	$$
\end{defn}

If (\ref{cl1}) is well-posed, we can close the loop with $\De_{ex}(\De)$ to get
$
\smat{\dot x_e\\z_p}
=\smat{\star&\star\\\star&0}
\smat{x_e\\w_p}
$
where the entries with $\star$ depend on $\De$ and $\De_c(\De)$; note that the structured LFRs (\ref{s00}), (\ref{s12}) imply (\ref{decp1}), (\ref{controller}) which lead to the desired zero block for $w_p\to z_p$ to render the $\Hz$-norm finite. Hence, the \textit{$\Hz$-gain-scheduling problem} involves a nontrivial structural requirement.
\begin{prob}\label{problem1}
	For a given bound $\ga>0$, determine a controller $\Kc(\De)$ structured as in (\ref{s12})-(\ref{dec}) such that
	\begin{enumerate}
		\item[(G1)] the controlled LFR (\ref{cl1}) is well-posed and stable,
		\item[(G2)] the squared $\Hz$-norm of $w_p\to z_p$ for linear time-varying systems (in the stochastic setting as in \cite{pagfer00}) is smaller than $\ga$ for $x_e(0)=0$ and for all $\De\in\mathbf{\De}$.
	\end{enumerate}
\end{prob}

\subsection{Analysis conditions for the original system}\label{sec22}

As well-known by the full block $S$-procedure, the conditions (G1)-(G2) are achieved if some matrix inequalities are feasible. This is formulated in the following standard analysis result from \cite{scherer2000} based on the class $\mathbf{\hat{P}}$ of full block scalings $\hat{\Pc}\hspace{0.1ex}{\in}\,\Sb^{(\uh+r^c+\vh+r^c)}$ satisfying
\begin{equation}\label{posconstr}
\sy \hat{\Pc}\smat{\De_{ex}(V)\\I_{\vh+r^c}}\cg0\quad \text{for all}\quad V\in\mathbf{V}.
\end{equation}

\begin{thm}\label{theo:analysisu}
	The design goals (G1)-(G2) are reached for the structured controller $K(\De)$ with (\ref{s12})-(\ref{dec}) if there exist $\Xc_1\cg0$, $Z\cg0$ with $\tr(Z)<1$ as well as $\hat{\Pc}\in\mathbf{\hat{P}}$ such that
	\begin{equation}\label{eq61}
	\begin{aligned}
	\Lc_{\text{sub}}&\left(\smat{-\Xc_1&0\\0&0},\hat{\Pc},P_Z,\smat{\hat{\Ac}_{ij}\\\hat{\Cc}_{j}}\right)\cl0,\\
	\Lc&\left(\smat{0&\Xc_1\\\Xc_1&0},\hat{\Pc},P_\ga,\smat{\hat{\Ac}_{ij}&\hat{\Bc}_{i}\\\hat{\Cc}_{j}&\hat{\Dc}}\right)\cl0
	\end{aligned}
	\end{equation}
	hold for the closed-loop system
	(\ref{cl1}) with
	\begin{equation}\label{pzga}
	P_Z:=\smat{0&0\\0&Z^{-1}},\quad P_\ga:=\smat{-\ga I&0\\0&0}.
	\end{equation}
\end{thm}
Since (\ref{eq61}) involve two inequalities with specific outer factors and $\hat{\Pc}$ is unstructured, we cannot directly eliminate or substitute the controller parameters for convexification. In the sequel, we thus introduce a novel design procedure, while, in view of \cite{scherer2000}, we anticipate the synthesis result to be formulated with the full block scaling class
\begin{equation}\label{primal}
\begin{aligned}
\mathbf{P}_p:=\left\{\right.&P\in\Sb^{\uh+\vh}\ \big|\
\sy P\smat{I_{\uh}\\0}\cl0\ \text{and}\\
&\left.\hspace{12ex} \sy P\smat{V\\I_{\vh}}\cg0\
\text{for all}\ V\in\mathbf{V} \right\}
\end{aligned}
\end{equation}
related to $\De$ and the corresponding dual scaling class
\begin{equation}\label{ddual}
\begin{aligned}
\mathbf{P}_d
:=\left\{\right.&\t{P}\in\Sb^{\uh+\vh}\ \big|\
\sy \t{P}\smat{0\\I_{\vh}}\cg0\ \text{and}\\
&\left.\hspace{8ex} \sy \t{P}\smat{I_{\uh}\\-V^T}\cl0\
\text{for all}\ V\in\mathbf{V} \right\}.
\end{aligned}
\end{equation}

\section{Lifting design procedure}\label{sec:lifting}

If $\mathbf{\hat{P}}$ is restricted in Theorem~\ref{theo:analysisu} to the class of positive real scalings $\smat{0&Q\\Q^T&0}$ satisfying the passivity condition related to (\ref{posconstr}), i.e. $\He\left[\de Q\right]\cg0$ for all real $\de\geq0$, the approach in \cite{roesinger2019b} shows that the anti-diagonal scaling block is a fundamental stumbling block for convexification by transformation.
This motivates to replace the intractable inequalities (\ref{eq61}) by a suitable, sufficient analysis condition for a certain class of passive scalings.

\subsection{Lifted plant and closed-loop formulation}\label{sec31}

First, let us define a new LFR by reformulating the equations for $G(\De)$ in (\ref{s00}). Note that $\hat{w}=\De\hat{z}$ is equivalent to $\hat{w}\,{=}\,-\hat{w}+2\De\hat{z}$ and thus to $w\,{=}\,\De_l(\De)z$ for $\De\,{\in}\,\mathbf{\De}$ where
\begin{equation}\label{delift}
w\,{:=}\,z\,{:=}\,\mat{c}{\hat{w}\\\hat{z}},\ \ \De_l(V)\,{:=}\,\mat{cc}{-I_{\uh}&2V\\0&I_{\vh}}\ \ \text{for}\ \ V\,{\in}\,\mathbf{V}.
\end{equation}
Similarly, we can rearrange the matrices in (\ref{s00}) related to the uncertainty channel $\hat{w}\to\hat{z}$ to infer that (\ref{s00}) is true iff
\begin{eqnarray}\label{s11}
\begin{aligned}
&\mat{c}{\dot{x}\\\hline z\\\hdl z_p\\ y}=\mat{c|c;{2pt/1pt}cc}{A_{11}&A_{12}&B_1^p&B_1\\\hline A_{21}&A_{22}&B_2^p&B_2\\\hdl
	C_1^p &C_2^p&D^p &D_1\\ C_1 & C_2& D_2&0}\mat{c}{x\\\hline w\\\hdl w_p\\ u}\\
&\hspace{-1ex}{:=}\,\mat{c|cc;{2pt/1pt}cc}{A_{11}&\Ah_{12}&0&B_1^p&B_1\\\hline 0&I_{\uh}&0&0&0\\ 2\Ah_{21}&2\Ah_{22}&-I_{\vh}&2\Bh_2^p&2\Bh_2\\\hdl
	C_1^p &\Ch_2^p&0&D^p &D_1\\ C_1 & \Ch_2& 0& D_2&0}\mat{c}{x\\\hline w\\\hdl w_p\\ u},\
w=\De_l(\De)z
\end{aligned}
\end{eqnarray}
holds for $\De\,{\in}\,\mathbf{\De}$. This construction results in a specifically structured uncertainty channel $w\,{\to}\,z$ of dimension $(\uh+\vh)\ti(\uh+\vh)$; in the sequel, we abbreviate (\ref{s11}) by $G_l(\De)$ and refer to $G_l(\De)/\De_l(\De)$ as \textit{lifted} LFR/\textit{lifted block}.\\
For the lifted LFR (\ref{s11}), let us describe the associated controller $K(\De)$ again by (\ref{s12})-(\ref{dec}) with the difference that $\De_c(.)$ is scheduled by $\De_l(\De)$ which, in general, leads to a larger size of the scheduling channel. For reasons of space, let use $\De_c(.)$ instead of $\De_c(\De_l(.))$ in the sequel.\\
By interconnecting (\ref{s11}) with (\ref{s12}), we get the closed-loop system
\begin{equation}\label{s13}
\mat{c}{\dot{x}_e\\\hdl z\\\hat{z}_c\\\hdl z_p}\mat{c;{2pt/1pt}c;{2pt/1pt}c}{\Ac_{11}&\Ac_{12}&\Bc_1\\\hdl \Ac_{21}&\Ac_{22}&\Bc_2\\\hdl \Cc_1&\Cc_2&\Dc}\mat{c}{x_e\\\hdl w\\\hat{w}_c\\\hdl w_p},\ \mat{c}{w\\\hat{w}_c}=\De_{lc}(\De)\mat{c}{z\\\hat{z}_c}
\end{equation}
with the corresponding scheduling block being defined as
\begin{equation}\label{delc}
\De_{lc}(V):=\mat{cc}{\De_l(V)&0\\0&\De_c(V)}\in\R^{(r^s+r^c)\ti (r^s+r^c)}
\end{equation}
for some $V\in\mathbf{V}$ and for the relevant dimensions
$$
n:=\nh+n^c,\ r^s:=\uh+\vh,\ r:=r^s+r^c=(\uh+\vh)+(r^c_1+r^c_2);
$$
the closed-loop matrices can be routinely expressed as
$$
\renewcommand{\arraystretch}{1.0}\arraycolsep=1.0pt
\scalebox{0.95}{$
\mat{c;{2pt/1pt}c}{\Ac_{ij}&\Bc_{i}\\\hdl \Cc_{j}&\Dc}\hspace{-0.3ex}=\hspace{-0.3ex}\mat{cc;{2pt/1pt}c}{A_{ij}&0&B^p_i\\0&0&0\\\hdl C^p_j&0&D^p}
+\mat{cc}{0&B_i\\I&0\\\hdl0&D_1}\hspace{-0.2ex}\mat{cc}{A^c_{ij}&B^c_i\\C^c_j&D^c}\hspace{-0.2ex}\mat{cc;{2pt/1pt}c}{0&I&0\\C_j&0&D_2}.$}
$$

\subsection{Lifted analysis conditions with passive scaling classes}

As a first observation, the scalings of $\mathbf{P}_p,\,\mathbf{P}_d$ in (\ref{primal})-(\ref{ddual}) already fulfill a passivity condition for the lifted block, i.e.
\begin{equation}\label{liftedprimal}
\begin{aligned}
&\hspace{-1.5ex}\mathbf{P}_p=\left\{P\in\Sb^{r^s}\, \big|\,
\He\left[P\De_l(V)\right]\cg0\ \, \text{for all}\ \, V\in\mathbf{V} \right\},\\
&\hspace{-1.5ex}\mathbf{P}_d=\left\{\t{P}\in\Sb^{r^s}\, \big|\,
\He\left[\t{P}\De_l(V)^T\right]\cg0\ \,  \text{for all}\ \, V\in\mathbf{V} \right\};
\end{aligned}
\end{equation}
this can be seen, e.g., for $\mathbf{P}_p$ by applying a congruence transformation with the invertible $\smat{I_{\uh}&V\\0&I_{\vh}}$ to the condition $\He\left[P\De_l(V)\right]\cg0$ for some $P\in\Sb^{\uh+\vh}$ and $V\in\mathbf{V}$. Secondly, if we replace $V$ by the lifted block $\De_l(V)$, the extended block $\De_{ex}(V)$ from Section~\ref{sec21} becomes $\De_{lc}(V)$ in (\ref{delc}). Hence, this motivates to define an appropriate scaling class $\mathbf{P}$ for the lifted $\De_{lc}(V)$ by some passivity condition as
\begin{equation}\label{scalingslifted}
\mathbf{P}:=\left\{\Pc\in\Sb^{r}\ \big|\ \He\left[\Pc\De_{lc}(V)\right]\cg0\hspace{1.85ex} \text{for all}\hspace{1.85ex} V\in\mathbf{V} \right\}.
\end{equation}
It will be insightful to see in Section~\ref{comp} that the specific choice of $\mathbf{P}$ causes no restriction if compared to the full block scaling class of \cite{scherer2000}.
Moreover, it will be crucial to see that a solution for Problem~\ref{problem1} can be obtained by solving the $\Hz$-gain-scheduling problem for the lifted LFR. This is achieved by starting, analogously to Section~\ref{sec22}, with the analysis inequalities
\begin{equation}\label{s14}
\begin{aligned}
\Lc_{\text{sub}}&\left(\smat{-\Xc_1&0\\0&0},\smat{0&\Pc\\\Pc&0},P_Z,\smat{\Ac_{ij}\\\Cc_{j}}\right)\cl0,\\
\Lc&\left(\smat{0&\Xc_1\\\Xc_1&0},\smat{0&\Pc\\\Pc&0},P_\ga,\smat{\Ac_{ij}&\Bc_{i}\\\Cc_{j}&\Dc}\right)\cl0
\end{aligned}
\end{equation}
for the controlled system (\ref{s13}) and scalings $\Pc\in\mathbf{P}$ with a passivity structure. As a crucial advantage over the original conditions for full block scalings in (\ref{eq61}), we show that (\ref{s14}) can be indeed convexified.

\subsection{Steps of lifting}\label{analysis}

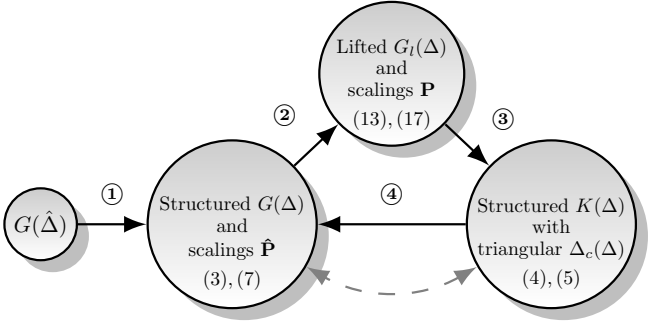
\begin{figure}
	\scalebox{2.1}{\hspace{-0.5ex}
	\begin{tikzpicture}
	\node[syr4] (g) at (0.8,0)  {};
	\node[syr7] (gs) at (2,0)  {};
	\node[syr6] (l) at (3,0.95)  {};
	\node[syr7] (c) at (4,0)  {};
	\node(gt) at (0.8,0) {\scalebox{0.4}{$G(\hat{\De})$}};
	\node(gst) at (2,0) {\scalebox{0.35}{$\begin{array}{c}\\[2ex]\text{Structured}\ G(\De)\\\text{and}\\\text{scalings}\ \mathbf{\hat{P}}\\[0.9ex](\ref{s00}),(\ref{posconstr})\end{array}$}};
	\node(lt) at (3,0.95) {\scalebox{0.35}{$\begin{array}{c}\\\text{Lifted}\ G_l(\De)\\\text{and}\\\text{scalings}\ \textbf{P}\\[0.9ex](\ref{s11}),(\ref{scalingslifted})\end{array}$}};
	\node(ct) at (4,0) {\scalebox{0.35}{$\begin{array}{c}\\[2ex]\text{Structured}\ K(\De)\\\text{with}\\ \text{triangular}\ \De_c(\De)\\[0.9ex](\ref{s12}),(\ref{dec})\end{array}$}};
	\draw[->] (g) -- node[pos=.5]{\scalebox{0.42}{\textcircled{\scriptsize\textbf{1}}}} (gs);
	\draw[->] (gs) -- node[pos=.5]{\scalebox{0.42}{\textcircled{\scriptsize\textbf{2}}}} (l);
	\draw[->] (l) -- node[pos=.6]{\scalebox{0.42}{\textcircled{\scriptsize\textbf{3}}}} (c);
	\draw[->] (c) -- node[pos=.5,swap]{\scalebox{0.42}{\textcircled{\scriptsize\textbf{4}}}} (gs);
	\draw [<->, bend angle=30, bend left,dashed,color=gray] (c) to (gs);
	\end{tikzpicture}}
	\caption{Steps of lifting technique: Build plant LFR $G(\De)$ in~\textcircled{\small\textbf{1}} and lifted plant LFR \textcircled{\small\textbf{2}}, design controller $K(\De)$ for lifted LFR \textcircled{\small\textbf{3}} and interconnect with $G(\De)$ in \textcircled{\small\textbf{4}}.}\label{blk2}
\end{figure}

Let us now summarize the concrete lifting design technique which is visualized in Fig.~\ref{blk2} and consists of four steps: For the first step \textcircled{\small\textbf{1}}, we have described in Section~\ref{sec1} the uncertain plant $G(\hat{\De})$ of (\ref{decp1}) by the structured LFR $G(\De)$ in (\ref{s00}) in order to formulate the analysis conditions (\ref{eq61}) with the class of full block scalings $\mathbf{\hat{P}}$. Next, we have performed the lifting step \textcircled{\small\textbf{2}} in Section~\ref{sec31} to  obtain the lifted LFR $G_l(\De)$ in (\ref{s11}). In the synthesis step \textcircled{\small\textbf{3}}, presented in Section~\ref{sec2}, we solve the associated $\Hz$-gain scheduling problem for the lifted LFR to obtain a structured controller $K(\De)$ with triangular $\De_c(\De)$ as in (\ref{s12})-(\ref{dec}). For this purpose, we rely on the analysis inequalities (\ref{s14})
for the lifted LFR and use the passive scaling class $\mathbf{P}$. The last step \textcircled{\small\textbf{4}} is given in Section~\ref{sec33} and clarifies that the constructed controller also solves the desired gain-scheduling Problem~\ref{problem1} for the original LFR $G(\De)$. Note that the design approach for positive real scalings in \cite{roesinger2019b} is only based on \textcircled{\small\textbf{1}} and the dashed grey lines, while \textcircled{\small\textbf{2}}-\textcircled{\small\textbf{4}} are the core novel synthesis steps for full block scalings.

\subsection{Consequences for the original system}\label{sec33}

The following result covers step \textcircled{\small\textbf{4}} of Fig.~\ref{blk2}.~
\begin{thm}\label{theo:analysis}
Suppose there exist a structured controller $K(\De)$ with (\ref{s12})-(\ref{dec}) as well as $\Xc_1\cg0$, $Z\cg0$ with $\tr(Z)<1$, $\Pc\in\mathbf{P}$ such that the closed-loop system (\ref{s13}) for the lifted LFR (\ref{s11}) fulfills (\ref{s14}) with $P_Z,\,P_\ga$ structured as in (\ref{pzga}).\\	
Then we can construct a full block scaling $\hat{\Pc}\in\mathbf{\hat{P}}$ with (\ref{posconstr}) such that the inequalities (\ref{eq61}) of Theorem~\ref{theo:analysisu} are true for the closed-loop system (\ref{cl1}) obtained for the initial plant LFR (\ref{s00}) and the same controller $K(\De)$.
\end{thm}

\begin{pf}
For some matrices $\Ahb$, $\Bhb$, $\Chb$, $Q\in\Sb$, $R\in\Sb$ and $S$ of suitable dimension, we first observe that
\begin{equation}\label{eqo}
\text{He}\left[\mat{c}{\Ahb\\\hdl\Chb}^T\hspace{-0.5ex}\mat{c;{2pt/1pt}c}{Q&S^T\\\hdl S&R}\mat{c}{\Ahb\\\hdl\Bhb}\right]
\,{=}\,\sy\mat{ccc}{2Q&S^T&S^T\\S&0&R\\S&R&0}\mat{c}{\Ahb\\\Bhb\\\Chb}.
\end{equation}
Now, let the analysis inequalities in (\ref{s14}) be satisfied for some $\Pc\in\mathbf{P}$ and for the lifted LFR interconnected with a given controller $K(\De)$. By the definition of $\mathbf{P}$, we infer $\He\left[\Pc\De_{lc}(V)\right]\cg0$ for all $V\in\mathbf{V}$. Applying for each $V\in\mathbf{V}$ a congruence transformation with	
$$
\mat{cc}{V&0\\I_{\vh}&0\\0&I_{r^c}}\ \ \text{yields}\ \
\He\left[\mat{cc}{V&0\\I_{\vh}&0\\\hdl 0&I_{r^c}}^T\hspace{-1ex}\Pc\mat{cc}{V&0\\I_{\vh}&0\\\hdl 0&\De_c(V)}\right]\cg0
$$
for all $V\in\mathbf{V}$. Next, let us partition $\Pc$ according to the outer factors of the latter inequality as
$$
\Pc=\mat{c;{2pt/1pt}c}{Q&S^T\\\hdl S&R}=\mat{cc;{2pt/1pt}c}{Q_{11}&Q_{12}&S_1^T\\Q_{21}&Q_{22}&S_2^T\\\hdl S_1&S_2&R}
$$
to conclude with (\ref{eqo}) after a suitable permutation that
\begin{equation}\label{eq9}
\sy
\underbrace{\mat{cc;{2pt/1pt}cc}{2Q_{11}&S_1^T&2Q_{12}&S_1^T\\S_1&0&S_2&R\\\hdl 2Q_{21}&S_2^T&2Q_{22}&S_2^T\\S_1&R&S_2&0}}_{=:\hat{\Pc}}
\mat{cc}{V&0\\0&\De_c(V)\\\hdl I_{\vh}&0\\0&I_{r^c}}\cg0.
\end{equation}
Thus $\hat{\Pc}\in\mathbf{\hat{P}}$. It is essential that the analysis inequalities (\ref{eq61}) obtained for (\ref{s00}) and for the same $K(\De)$ are also valid for the constructed $\hat{\Pc}$ from (\ref{eq9}). This follows by applying suitable congruence transformations to (\ref{s14}) along with (\ref{eqo}); we need to omit the details for reasons of space.\hfill$\blacksquare$
\end{pf}

\subsection{Comparison of scaling classes}\label{comp}

Let $\mathbf{\hat{P}_F}$ be the full block scaling class used for gain-scheduling in \cite{scherer2000}. Note that $\mathbf{\hat{P}_F}$ is a subset of $\mathbf{\hat{P}}$ from Section~\ref{sec22} and consists of all scalings $\hat{\Pc}\in\Sb^{(\uh+r^c+\vh+r^c)}$ satisfying in addition to (\ref{posconstr}) the constraints
\begin{equation}\label{addconstraints}
\sy\hat{\Pc}\smat{I_{\uh+r^c}\\0}\cl0\quad\text{and}\quad  \sy\hat{\Pc}\smat{0\\I_{\vh+r^c}}\cg0.
\end{equation}
We emphasize that it is not at all clear how to convexify the synthesis problem based on (\ref{eq61}) for the class $\mathbf{\hat{P}_F}$. Still, let us briefly sketch that the choice of the specifically structured scalings $\mathbf{P}$ in (\ref{scalingslifted}) causes no conservatism, i.e., if $\ga_F$ is the optimal bound obtained for (\ref{eq61}) with $\mathbf{\hat{P}_F}$, and $\ga_l$ denotes the one for synthesis based on (\ref{s14}) with the lifted LFR and $\mathbf{P}$, the relation $\ga_l\leq\ga_F$ always holds.\\
For this purpose, let us perform the lifting step in Section~\ref{sec31} both for the plant $G(\De)$ and for $K(\De)$.
This leads to the lifted plant LFR $G_l(\De)$ in (\ref{s11}) as well as to a lifted controller LFR $K_l(\De)$ with a scheduling channel resembling the structure of those for $G_l(\De)$, while being scheduled by the structured $\De_l(\De_c(\De))$ with $\De_l(.)$ from (\ref{delift}). Note that the resulting LFR of $K_l(\De)$ can be always obtained by a structural restriction of the LFR matrices for $K(\De)$.
By exploiting the scaling properties (\ref{addconstraints}), (\ref{posconstr}) imposed for $\mathbf{\hat{P}_F}$, it is crucial to see that the original analysis inequalities (\ref{eq61}) hold for some $\hat{\Pc}\in\mathbf{\hat{P}_F}$ if and only if the modified analysis inequalities (\ref{s14}) are satisfied for the closed-loop system obtained from interconnecting $G_l(\De)$ with the lifted controller LFR $K_l(\De)$, and for some scaling $\Pc\in\Sb$ satisfying the passivity constraint
\begin{equation}\label{addconstraint2}
\text{He}\left[\Pc\smat{\De_l(V)&0\\0&\De_l(\De_c(V))}\right]\cg0\quad\text{for all}\quad V\in\mathbf{V}.
\end{equation}
We omit the details for reasons of space, but remark that, upon permutation, $\Pc$ in (\ref{addconstraint2}) equals $\hat{\Pc}$.
We observe that (\ref{addconstraint2}) is exactly the condition that appears for the passive scalings $\mathbf{P}$ in (\ref{scalingslifted}) if replacing $\De_l(\De_c(V))$ by $\De_c(V)$. Since the class of all LFRs for $K(\De)$ encompasses that of all LFRs for $K_l(\De)$ as argued above, we infer $\ga_l\leq\ga_F$.

\section{Synthesis for lifted system}\label{sec2}

In the following part we deal with the synthesis step \textcircled{\small\textbf{3}} in Fig.~\ref{blk2}, i.e., we use a structured controller parameter transformation combined with a suitable scaling factorization to solve the $\Hz$-gain-scheduling problem for the lifted LFR. In the context of structured $\Hz$-design, a related factorization is established for positive definite matrices in \cite{scherer2014} to design triangular, time-invariant controllers, as well as for positive real matrices in \cite{roesinger2019b} to synthesize gain-scheduled controllers with a diagonal scheduling function of scalar parameters.
Technically, we show as a novel step that the passivity condition for $\mathbf{P}$ in (\ref{scalingslifted}) can be used to derive a structured factorization for possibly indefinite scalings which is used to guarantee the existence of a block-triangular scheduling function $\De_c(.)$ for matrix parameters.\\
Before formulating the main synthesis result, we present the corresponding variables which consist of the matrices
$X_1, Y_1\in\Sb^{\nh}$.
Further we take
\begin{equation}\label{svar1}
X_2=\mat{cc}{Q_2&Q_3}\qquad\text{and}\qquad Y_2=\mat{cc}{\Qa_1&I_{r^s}}
\end{equation}
of dimension $r^s\ti(r^s+r^s)$ with $Q_2$, $Q_3\in\mathbf{P}_p$, $\Qa_1\in\mathbf{P}_d$ of dimension $r^s\ti r^s$ where the sets $\mathbf{P}_p$ and $\mathbf{P}_d$ are given in (\ref{liftedprimal}).
Moreover, for a compact notation, we use
\begin{equation}\label{rec}
\renewcommand{\arraystretch}{1.2}
\mat{c;{2pt/1pt}c;{2pt/1pt}c}{\KB_{11}&\KB_{12}&\LB_1\\\hdl \KB_{21}&\KB_{22}&\LB_2\\\hdl \MB_1&\MB_2&\NB}
:=\mat{c;{2pt/1pt}cc;{2pt/1pt}c}{\Ks_{11}&\Ks_{12}&\Ks_{13}&\Ls_1\\\hdl \Ks_{21}&\Ks_{22}&Q_2^TA_{22}&0\\
	\Ks_{31}&\Ks_{32}&\Ks_{33}&\Ls_3\\\hdl \Ms_1&\Ms_2&0&0}
\end{equation}
of the partition $(\nh+(r^s+r^s)+m)\ti (\nh+(r^s+r^s)+k)$ which includes $Q_2$ from (\ref{svar1}) and the unstructured variables $\Ks_{ij}$, $\Ls_i$, $\Ms_j$. This leads to the following $\Hz$-gain scheduling synthesis result.
\begin{thm}\label{theo1}
	Let $\ga>0$ be fixed. There exists a structured controller with triangular scheduling function $\De_c(.)$ as in (\ref{s12})-(\ref{dec}) and some $\Xc_1\cg0$, $\Pc\in\mathbf{P}$, $Z\cg0$ with $\tr(Z)<1$ such that the inequalities (\ref{s14}) (with (\ref{pzga})) hold for the closed-loop system (\ref{s13}) iff there exist $X_1,\,Y_1\in\Sb^{\nh}$, structured $X_2$, $Y_2$ from (\ref{svar1}), $\Ks_{ij},\,\Ls_i,\,\Ms_j$ with (\ref{rec}), and some $Z\cg0$ with $\tr(Z)<1$ such that
	\begin{equation}\label{s211}
	\begin{aligned}
	\Lc_{\text{sub}}&\left(\smat{-\XB&0\\0&0},\smat{0&I\\I&0},P_Z,\smat{\AB_{ij}\\\CB_{j}}\right)\cl0,\\
	\Lc&\left(\smat{0&I\\I&0},\smat{0&I\\I&0},P_\ga,\smat{\AB_{ij}&\BB_{i}\\\CB_{j}&\DB}\right)\cl0
	\end{aligned}
	\end{equation}
	are satisfied after inserting for $i,j=1,2$ the blocks
	\begin{equation}\label{s22}
	\begin{aligned}
	&\hspace{9.55ex}\XB:=\mat{cc}{Y_1&I_{\nh}\\I_{\nh}&X_1},\\
	&\mat{c;{2pt/1pt}c}{\AB_{ij}&\BB_i\\\hdl \CB_j&\DB}
	:=\mat{cc;{2pt/1pt}c}{A_{ij}Y_j&A_{ij}&B_i^p\\0&X_i^TA_{ij}&X_i^TB_i^p\\\hdl C_j^pY_j&C_j^p&D^p}+\\&\hspace{13ex}+
	\mat{cc}{0&B_i\\I&0\\\hdl 0&D_1}
	\mat{cc}{\KB_{ij}&\LB_i\\\MB_j&\NB}
	\mat{cc;{2pt/1pt}c}{I&0&0\\0&C_j&D_2}.
	\end{aligned}
	\end{equation}
\end{thm}
Since $\mathbf{V}=\text{Co}\{\De_1,\ldots,\De_N\}$ and the sets $\mathbf{P}_p$, $\mathbf{P}_d$ can be expressed as in (\ref{primal}), (\ref{ddual}), the conditions $Q\in\mathbf{P}_p$, $\t{Q}\in\mathbf{P}_d$ reduce to finitely many inequalities (see \cite{scherer2000}):
$$
\begin{aligned}
\sy Q\smat{I\\0}&\cl0, &\ \ \sy Q\smat{\De_i\\I}&\cg0,\\
\sy\t{Q}\smat{0\\I}&\cg0, &\ \ \sy\t{Q}\smat{I\\-\De_i^T}&\cl0\quad\text{for}\quad i=1,\ldots,N.
\end{aligned}
$$
After applying the Schur complement to (\ref{s211}), we get a standard LMI test with finitely many constraints such that a direct minimization over $\gamma$ is possible.
We present the proof of Theorem~\ref{theo1} in Appendix~\ref{sproof}. Note that our proof is constructive, i.e., if the associated LMIs are feasible, a suitable $\Hz$-controller (\ref{s12})-(\ref{dec}) can be constructed with McMillan degree of at most $\nh$ and scheduling block size $r^c$ of at most $2r^s$, while we give an explicit formula for~$\De_c(.)$.

\begin{rem}
Analogously to Remark~5 and 6 in \cite{roesinger2019b}, Theorem~\ref{theo1} can also handle gain-scheduling with quadratic performance and multiple objectives by properly modifying $P_\ga$. Also $\Ks_{11},\Ks_{12},\Ks_{13},\Ls_1$ can be partially eliminated to reduce the number of variables.
\end{rem}

\section{A numerical example}\label{sec3}

\begin{figure}
	\hspace{-3.8ex}
	\includegraphics[width=0.53\textwidth]{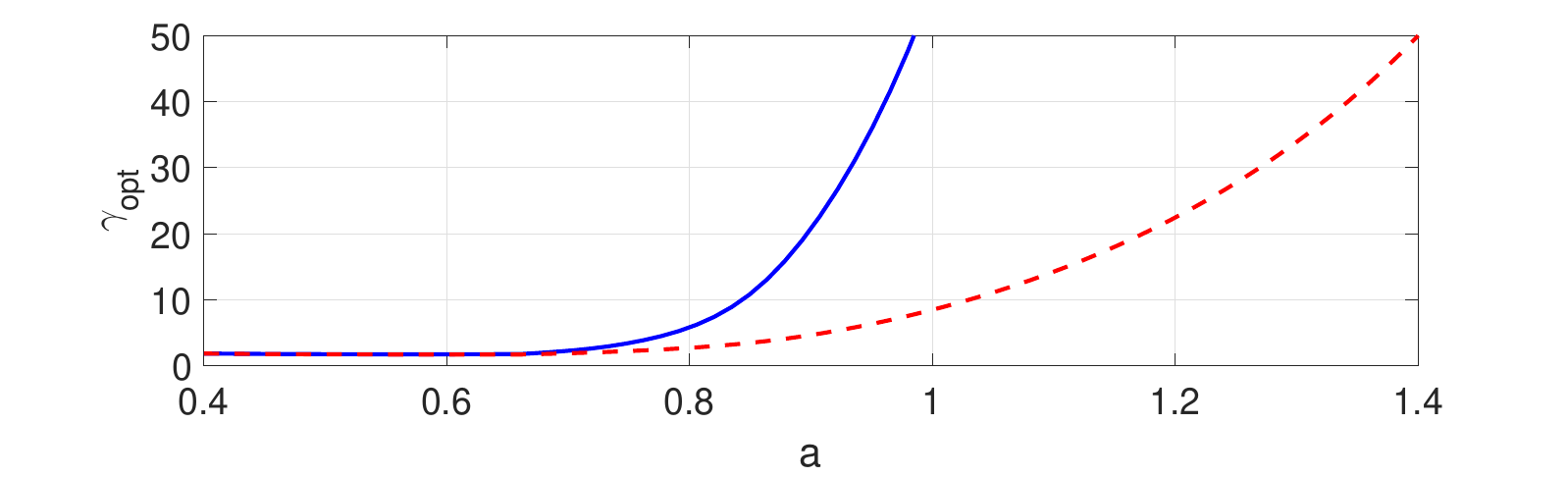}
	\caption{Optimal bounds $\ga_{\text{opt}}$ for the lifted design (dashed red) and $D/G$-scalings (full blue) with $a\in[0.4,1.4]$.}
	\label{plot1}
\end{figure}

To present a short academic example, let the matrices of the structured LFR in (\ref{s00}) be given as in Section~4.2 of \cite{roesinger2019b} with $\hat{A}_{12}$ depending on some parameter $a\in[0.4,1.4]$. Moreover, let $\De=\diag(\de_1I_2,\de_2)$ be of size $3\ti3$ with time-varying parametric uncertainties $\de_1(t)\in[-0.8,0.8]$, $\de_2(t)\in[-0.6,0.6]$. Based on implementations of our algorithms in the Matlab Robust Control Toolbox, we compare in Fig.~\ref{plot1} the optimal bounds $\ga_{\text{opt}}$ of the squared $\Hz$-norm for the lifted design (dashed red) obtained for the passive scaling class $\mathbf{P}$ from (\ref{scalingslifted}) with $D/G$-scalings (full blue). Note that $\Hz$-gain-scheduling synthesis for $D/G$-scalings with structured LFRs can be performed with the positive real scaling results from \cite{roesinger2019b} for the original LFR (\ref{s00}) by using the well-known M\"obius transformation to map the uncertainty intervals for $\de_i$ into $[0,\infty]$. To the best knowledge of the authors, there exist no alternative approaches that solve the underlying structured $\Hz$-design problem in this generality. The results confirm that the lifted approach is less conservative than $D/G$-scalings as expected from Section~\ref{comp}. In particular, beyond the shown parameter range for $a$, the synthesis LMIs get infeasible for $D/G$-scalings if $a$ approaches $1.67$, while the lifted design is feasible up to $a=2.17$.

\section{Conclusion and outlook}

In this work, we have introduced a new lifting technique to synthesize controllers for the $\Hz$-gain-scheduling problem with full block scalings. Especially, our design framework guarantees finiteness of the closed-loop $\Hz$-norm by relying on structured plant and controller LFRs, and by constructing a block-triangular scheduling function. We hope that these new methodologies offer manifold potential for refined synthesis results as the combination with parameter-dependent Lyapunov functions. A further task is the investigation of possible numerical advantages of the used scaling extension over existing approaches.

\bibliography{lit}

\begin{thebibliography}{21}
\providecommand{\natexlab}[1]{#1}
\providecommand{\url}[1]{\texttt{#1}}
\providecommand{\urlprefix}{URL }
\expandafter\ifx\csname urlstyle\endcsname\relax
  \providecommand{\doi}[1]{doi:\discretionary{}{}{}#1}\else
  \providecommand{\doi}{doi:\discretionary{}{}{}\begingroup
  \urlstyle{rm}\Url}\fi

\bibitem[{Apkarian and Adams(1997)}]{apkarian1997}
Apkarian, P. and Adams, R.J. (1997).
\newblock Advanced gain-scheduling techniques for uncertain systems.
\newblock In \emph{Proc. American Control Conf.}, 3331--3335.

\bibitem[{Apkarian and Gahinet(1995)}]{apkarian1995}
Apkarian, P. and Gahinet, P. (1995).
\newblock A convex characterization of gain-scheduled $\mathcal{H}_\infty$
  controllers.
\newblock \emph{IEEE Trans. Autom. Control}, 40(5), 853--864.

\bibitem[{Becker(1995)}]{becker1995}
Becker, G. (1995).
\newblock Parameter-dependent control of an under-actuated mechanical system.
\newblock In \emph{Proc. 34th IEEE Conf. Decision and Control}, 543--548.

\bibitem[{de~Souza and Trofino(2006)}]{souza2006}
de~Souza, C.E. and Trofino, A. (2006).
\newblock Gain-scheduled $\mathcal{H}_2$ controller synthesis for linear
  parameter varying systems via parameter-dependent {Lyapunov} functions.
\newblock \emph{Int. J. Robust Nonlin.}, 16(5), 243--257.

\bibitem[{Guerreiro et~al.(2007)Guerreiro, Silvestre, Cunha, and
  Antunes}]{guerreiro2007}
Guerreiro, B., Silvestre, C., Cunha, R., and Antunes, D. (2007).
\newblock Trajectory tracking $\mathcal{H}_2$ controller for autonomous
  helicopters: {An} application to industrial chimney inspection.
\newblock \emph{IFAC Proc. Vol.}, 40(7), 431--436.

\bibitem[{Helmersson(1998)}]{helmersson98}
Helmersson, A. (1998).
\newblock $\mu$ synthesis and {LFT} gain scheduling with real uncertainties.
\newblock \emph{Int. J. Robust Nonlin.}, 8(7), 631--642.

\bibitem[{Langbort et~al.(2004)Langbort, Chandra, and D'Andrea}]{langbort2004}
Langbort, C., Chandra, R.S., and D'Andrea, R. (2004).
\newblock Distributed control design for systems interconnected over an
  arbitrary graph.
\newblock \emph{IEEE Trans. Autom. Control}, 49(9), 1502--1519.

\bibitem[{Masubuchi et~al.(1998)Masubuchi, Ohara, and Suda}]{masubuchi1998}
Masubuchi, I., Ohara, A., and Suda, N. (1998).
\newblock {LMI}-based controller synthesis: A unified formulation and solution.
\newblock \emph{Int. J. Robust Nonlin.}, 8(8), 669--686.

\bibitem[{Mustaki et~al.(2019)Mustaki, Nguyen, Chevrel, Yagoubi, and
  Fauvel}]{mustaki2019}
Mustaki, S., Nguyen, A.T., Chevrel, P., Yagoubi, M., and Fauvel, F. (2019).
\newblock Comparison of two robust static output feedback {$H_2$} design
  approaches for car lateral control.
\newblock In \emph{Proc. 18th European Control Conf.}, 716--723.

\bibitem[{Packard(1994)}]{packard1994}
Packard, A. (1994).
\newblock Gain-scheduling via linear fractional transformations.
\newblock \emph{Syst. Control Lett.}, 22(2), 79--92.

\bibitem[{Paganini and Feron(2000)}]{pagfer00}
Paganini, F. and Feron, E. (2000).
\newblock Linear matrix inequality methods for robust $\mathcal{H}_2$ analysis:
  {A} survey with comparisons.
\newblock In L.~El~Ghaoui and S.I. Niculescu (eds.), \emph{Advances in linear
  matrix inequality methods in control}, 129--151. SIAM, Philadelphia.

\bibitem[{R\"osinger and Scherer(2019)}]{roesinger2019b}
R\"osinger, C.A. and Scherer, C.W. (2019).
\newblock A scalings approach to $\mathcal{H}_2$-gain-scheduling synthesis
  without elimination.
\newblock \emph{IFAC-PapersOnLine}, 52(28), 50--57.

\bibitem[{Sato(2011)}]{sato2011}
Sato, M. (2011).
\newblock Gain-scheduled output-feedback controllers depending solely on
  scheduling parameters via parameter-dependent {L}yapunov functions.
\newblock \emph{Automatica}, 47(12), 2786--2790.

\bibitem[{Scherer(2000)}]{scherer2000}
Scherer, C.W. (2000).
\newblock Robust mixed control and linear parameter-varying control with full
  block scalings.
\newblock In L.~El~Ghaoui and S.I. Niculescu (eds.), \emph{Advances in linear
  matrix inequality methods in control}, 187--207. SIAM, Philadelphia.

\bibitem[{Scherer(2014)}]{scherer2014}
Scherer, C.W. (2014).
\newblock {$H_\infty$}- and {$H_2$}-synthesis for nested interconnections: A
  direct state-space approach by linear matrix inequalities.
\newblock In \emph{21st Int. Symp. Math. Theory of Netw. and Syst.},
  1589--1594.

\bibitem[{Scherer et~al.(1997)Scherer, Gahinet, and Chilali}]{scherer1997}
Scherer, C.W., Gahinet, P., and Chilali, M. (1997).
\newblock Multi-objective output-feedback control via {LMI} optimization.
\newblock \emph{IEEE Trans. Autom. Control}, 42(7), 896--911.

\bibitem[{Scorletti and {El~Ghaoui}(1998)}]{scorletti1998}
Scorletti, G. and {El~Ghaoui}, L. (1998).
\newblock {Improved LMI conditions for gain scheduling and related control
  problems}.
\newblock \emph{Int. J. Robust Nonlin.}, 8(10), 845--877.

\bibitem[{Veenman and Scherer(2013)}]{veenman2013}
Veenman, J. and Scherer, C.W. (2013).
\newblock Stability analysis with integral quadratic constraints: A
  dissipativity based proof.
\newblock In \emph{Proc. 52nd IEEE Conf. Decision and Control}, 3770--3775.

\bibitem[{Veenman and Scherer(2014)}]{veenman2014}
Veenman, J. and Scherer, C.W. (2014).
\newblock A synthesis framework for robust gain-scheduling controllers.
\newblock \emph{Automatica}, 50(11), 2799--2812.

\bibitem[{Wu and Dong(2005)}]{wu2005}
Wu, F. and Dong, K. (2005).
\newblock {Gain-scheduling control of LFT systems using parameter-dependent
  Lyapunov functions}.
\newblock In \emph{Proc. American Control Conf.}, 587--592.

\bibitem[{Wu et~al.(1996)Wu, Yang, Packard, and Becker}]{wu1996}
Wu, F., Yang, X.H., Packard, A., and Becker, G. (1996).
\newblock Induced {$L_2$}-norm control for {LPV} systems with bounded parameter
  variation rates.
\newblock \emph{Int. J. Robust Nonlin.}, 6(9-10), 983--998.

\end{thebibliography}

\appendix

\section{Proof of Theorem~\ref{theo1}}\label{sproof}

\textit{Necessity.}
	Let (\ref{s14}) be satisfied for (\ref{s13}), $\Xc_1\cg0$, $Z\cg0$ with $\tr(Z)<1$, and $\Xc_2:=\Pc\in\mathbf{P}$, i.e.
		\begin{equation}\label{p0}
		\begin{aligned}
		\Lc_{\text{sub}}&\left(\smat{-\Xc_1&0\\0&0},\smat{0&I\\I&0},P_Z,\smat{\Xc_i\Ac_{ij}\\\Cc_{j}}\right)\cl0,\\
		\Lc&\left(\smat{0&I\\I&0},\smat{0&I\\I&0},P_\ga,\smat{\Xc_i\Ac_{ij}&\Xc_i\Bc_{i}\\\Cc_{j}&\Dc}\right)\cl0.
		\end{aligned}
		\end{equation}
	\textbf{Step 1} (\textit{Factorizations}).\\
	W.l.o.g., let us assume that $n^c\geq \nh$ to factorize $\Xc_1$ as
	\begin{equation}\label{fac}
	\Xc_i\Yc_i=\Zc_i\ \ \text{with}\ \ \Yc_i:=\mat{cc}{Y_i&I\\V_i&0},\ \Zc_i:=\mat{cc}{I&X_i\\0&U_i}
	\end{equation}
	for $i=1$ such that $\Yc_1$ has full column rank (see \cite{scherer1997}).\\
	Moreover, if we assume that $r^c_1\geq r^s$ and $r^c_2\geq r^s$, let us show that $\Xc_2$ can be also factorized as in (\ref{fac}) such that $\Yc_2$ has full column rank where $V_2$ and $U_2$ are lower and upper block-triangular matrices, respectively, with respect to the partition $(r^c_1+r^c_2)\ti(r^s+r^s)$, and where $X_2$, $Y_2$ are partitioned as in (\ref{svar1}) for some suitable blocks $Q_2$, $Q_3$, $\Qa_1$.\\
	For this purpose, let us first clarify that
	$\Xc_2\in\mathbf{P}$ is invertible with some sub-blocks of full column rank, while we use the following partitions according to $r=r^s+r^c_1+r^c_2$:
	\begin{equation}\label{p1}
	\scalebox{0.95}{$
	\renewcommand{\arraystretch}{1.25}
	\Xc_2=\mat{c;{2pt/1pt}cc}{Q_3&S_{13}^T&S_{23}^T\\\hdl S_{13}&R_{11}&R_{21}^T\\S_{23}&R_{21}&R_{22}},\
	\Xc_2^{-1}\,{=}\,\mat{c;{2pt/1pt}cc}{\Qa_{1}&\Sa_{11}^T&\Sa_{21}^T\\\hdl \Sa_{11}&\Ra_{11}&\Ra_{21}^T\\\Sa_{21}&\Ra_{21}&\Ra_{22}}.
	$}
	\end{equation}
	For the given partition of $\Xc_2$ in (\ref{p1}), we note that $S_{13}$, $S_{23}$ are tall due to $r^c_j\geq r^s$ for $j=1,2$.
	Let us firstly perturb $R_{11},\,R_{21},\,R_{22}$ to achieve invertibility of $R_{22}$ and $\smat{R_{11}&R_{21}^T\\R_{21}&R_{22}}$. This allows to perturb $S_{13},\,S_{23},\,Q_3$ such that
	\begin{equation}\label{HSa}
	H:=-\mat{cc}{I&0}\mat{cc}{R_{11}&R_{21}^T\\R_{21}&R_{22}}^{-1}\mat{c}{S_{13}\\S_{23}},\ \Sa_{22}:=-R_{22}^{-1}S_{23}
	\end{equation}
	have full column rank and $Q_3-\smat{S_{13}^T&S_{23}^T}
	\smat{R_{11}&R_{21}^T\\R_{21}&R_{22}}^{-1}\smat{S_{13}\\S_{23}}$ is invertible.
	In particular, this implies invertibility of $\Xc_2$.
	Immediately, we infer that (\ref{fac}) is true for $i=2$ with
	\begin{equation}\label{p2}
	\scalebox{0.95}{$
	\renewcommand{\arraystretch}{1.25}
	\mat{c}{X_2\\\hdl U_2}\,{:=}\,\mat{cc}{Q_2&Q_3\\\hdl S_{12}&S_{13}\\0&S_{23}},\quad
	\mat{c}{Y_2\\\hdl V_2}\,{:=}\,\mat{cc}{\Qa_1&I_{r^s}\\\hdl\Sa_{11}&0\\\Sa_{21}&\Sa_{22}}$}
	\end{equation}
	where $Q_2:=Q_3-S_{23}^TR_{22}^{-1}S_{23}$ and $S_{12}:=S_{13}-R_{21}^TR_{22}^{-1}S_{23}$.
	By the block-inversion formula, we note that $\Qa_{1}$ is inver\-tible which, combined with (\ref{HSa}), reveals that
	$\Sa_{11}=H\Qa_{1}$ and $\Sa_{22}$ have full column rank. Thus, $V_2$ has full column rank which implies the same for $\Yc_2$ in (\ref{fac}).\\
	\textbf{Step 2} (\textit{Proof that $Q_2$, $Q_3\in\mathbf{P}_p$, $\Qa_1\in\mathbf{P}_d$}).\\
	For brevity, let us omit the argument of $\De_l(.),\,\De_c(.)$ and $\De_{lc}(.)$. Further, let us split $\De_{lc}$ into two parts such that
	\begin{equation}\label{p31}
	\begin{aligned}\renewcommand{\arraystretch}{0.9}\arraycolsep=1.0pt
	0&\cl \He\left[\Xc_2\De_{lc}\right]=\He\left[\Xc_2\mat{cc}{\De_l&0\\ 0&0}\right]+\He\left[\Xc_2\mat{cc}{0&0\\0&\De_c}
	\right].
	\end{aligned}
	\end{equation}
	Let us perform a congruence transformation with $\Yc_2$ on (\ref{p31}) while using (\ref{fac}) for $i=2$ and (\ref{p2}). This leads to 
	\begin{equation}\label{cond_scheduling}
	\begin{aligned}
	&\hspace{-1.47ex}0\cl\renewcommand{\arraystretch}{1.25}
	\mat{cc;{2pt/1pt}c}{\He\left[\Qa_1\De_l^T\right]&\ast&\ast\\\hdl Q_2\De_l\Qa_1+\De_l^T&\He\left[Q_2\De_l\right]&\ast\\ Q_3\De_l\Qa_1+\De_l^T&Q_3\De_l+\De_l^TQ_2&\He\left[Q_3\De_l\right]}
	+\\
	&\hspace{2ex}+\He\left[\renewcommand{\arraystretch}{1.25}\mat{c}{0\\\hdl U_2^T}\De_c
	\mat{c;{2pt/1pt}c}{V_2&0}\right].
	\end{aligned}
	\end{equation}
	Since $U_2^T,\,V_2,\,\De_c$ are lower block-triangular, the diagonal entries of (\ref{cond_scheduling}) just read as $Q_2$, $Q_3\in\mathbf{P}_p$, $\Qa_1\in\mathbf{P}_d$.\\
	\textbf{Step 3} (\textit{Derivation of synthesis inequalities (\ref{s211})}).\\
	Let us use the factorizations in (\ref{fac}) to apply congruence transformations with $\Yc_i$ to (\ref{p0}) for $i=1,2$. We get
		\begin{equation}\label{s17}
		\begin{aligned}
		\Lc_{\text{sub}}&\left(\smat{-\Zc_1^T\Yc_1&0\\0&0},\smat{0&I\\I&0},P_Z,\smat{\Zc_i^T\Ac_{ij}\Yc_j\\\Cc_{j}\Yc_j}\right)\,{\cl}\,0,\\	
		\Lc&\left(\smat{0&I\\I&0},\smat{0&I\\I&0},P_\ga,\smat{\Zc_i^T\Ac_{ij}\Yc_j&\Zc_i^T\Bc_{i}\\\Cc_{j}\Yc_j&\Dc}\right)\,{\cl}\,0.
		\end{aligned}
		\end{equation}
	By matching (\ref{s17}) to (\ref{s211}), the necessity part can then be finished similarly to \cite{roesinger2019b}: By symmetry, $\Zc_1^T\Yc_1$ equals $\XB$ from (\ref{s22}). Further, some calculations reveal that
	$$
	\renewcommand{\arraystretch}{0.9}\arraycolsep=2.0pt
	\begin{aligned}
	\mat{c;{2pt/1pt}c}{\Zc^T_i\Ac_{ij}\Yc_j&\Zc_i^T\Bc_i\\\hdl \Cc_j\Yc_j&\Dc}
	&=\mat{cc;{2pt/1pt}c}{A_{ij}Y_j&A_{ij}&B_i^p\\0&X_i^TA_{ij}&X_i^TB_i^p\\\hdl C_j^pY_j&C_j^p&D^p}+\\&\hspace{2.6ex}+
	\mat{cc}{0&B_i\\I&0\\\hdl 0&D_1}
	\mat{cc}{K_{ij}&L_i\\M_j&N}
	\mat{cc;{2pt/1pt}c}{I&0&0\\0&C_j&D_2}
	\end{aligned}
	$$
	for $i,j=1,2$ after performing the substitution
	\begin{equation}\label{s2}
	\scalebox{0.95}{$
	\begin{aligned}
	&\hspace{-4ex}\mat{c;{2pt/1pt}c}{K_{ij}&L_i\\\hdl M_j&N}:=
	\mat{c;{2pt/1pt}c}{X_i^TA_{ij}Y_j&0\\\hdl 0&0}+\\&\hspace{7ex}+
	\mat{cc}{U_i^T&X_i^TB_i\\\hdl 0&I_m}
	\mat{cc}{A_{ij}^c&B^c_i\\ C^c_j&D^c}
	\mat{c;{2pt/1pt}c}{V_j&0\\ C_jY_j&I_k}.
	\end{aligned}$}
	\end{equation}
	Moreover, by exploiting the sparsity structure of the controller matrices and $U_2$, $V_2$, we can introduce
	\begin{equation}\label{s3}
	\hspace{-0.9ex}
	\scalebox{0.93}{$
	\renewcommand{\arraystretch}{1.25}\arraycolsep=1.0pt
	\mat{c;{2pt/1pt}cc;{2pt/1pt}c}{\Ks_{11}&\Ks_{12}&\Ks_{13}&\Ls_1\\\hdl \Ks_{21}&\Ks_{22}&0&0\\\Ks_{31}&\Ks_{32}&\Ks_{33}&\Ls_3\\\hdl\Ms_1&\Ms_2&0&0}\,{:=}\,
	\mat{c;{2pt/1pt}c;{2pt/1pt}c}{K_{11}&K_{12}&L_1\\\hdl K_{21}&K_{22}-\smat{0&Q_2^TA_{22}\\0&0}&L_2\\\hdl M_1&M_2&0}
	$}
	\end{equation}
	which shows that (\ref{s17}) can be rewritten as (\ref{s211}) for (\ref{s22}).
	
	\textit{Sufficiency.}
	Let the inequalities in (\ref{s211}) be satisfied for (\ref{s22}) which comprises $X_1,\,Y_1\in\Sb^{\nh}$, structured $X_2$, $Y_2$ from (\ref{svar1}), $\Ks_{ij},\,\Ls_i,\,\Ms_j$ with (\ref{rec}), and $Z\cg0$ with $\tr(Z)<1$.\\
	\textbf{Step 1} (\textit{Construction of $\Xc_1$ and $\Xc_2$}).\\
	To define $\Yc_1,\,\Zc_1$ by (\ref{fac}), we choose $U_1:=I_{\nh}$, $V_1:=I_{\nh}-X_1^TY_1$. Hence, $\Zc_1^T\Yc_1=\XB$ and, since $\XB\cg0$ holds by (\ref{s211}), the matrices $U_1$, $V_1$ are invertible which implies the same for $\Yc_1,\,\Zc_1$. Thus (\ref{fac}) holds for $i=1$ with $\Xc_1:=\Zc_1\Yc_1^{-1}$.\\
	To find some suitable $\Xc_2$, we can achieve invertibility of
	$$
	T_1:=Q_2-\Qa_1^{-1}\quad \text{and}\quad T_2:=Q_3-Q_2
	$$
	by perturbation. For any invertible matrices $S_{13}$, $S_{23}$, let
	$$
	R_{21}=0\quad \text{and}\quad R_{11}:=S_{13}T_1^{-1}S_{13}^T,\quad R_{22}:=S_{23}T_2^{-1}S_{23}^T.
	$$
	This shows the validity of $Q_2=Q_3-S_{23}^TR_{22}^{-1}S_{23}$ and
	$$
	Q_3-\mat{cc}{S_{13}^T&S_{23}^T}\mat{cc}{R_{11}&R_{21}^T\\R_{21}&R_{22}}^{-1}\mat{c}{S_{13}\\S_{23}}\\
	=\Qa_{1}^{-1}.
	$$
	Therefore, if we define $\Xc_2$ by the first relation in (\ref{p1}), the block-inversion formula reveals that $\Xc_2$ is invertible with its inverse satisfying the second relation in (\ref{p1}) for some suitable $\Sa_{11},\,\Sa_{21},\,\Ra_{11},\,\Ra_{21},\,\Ra_{22}$. Further, let us take $\Sa_{22}:=-R_{22}^{-1}S_{23}$, $S_{12}:=S_{13}$ to define $U_2,\,V_2$ by (\ref{p2}) and $\Yc_2,\,\Zc_2$ by (\ref{fac}). Hence, (\ref{fac}) is true for $i=2$. Moreover, we identify $\Sa_{11}=-S_{13}^{-T}T_1\Qa_{1}$ which shows that $U_2,\,V_2$ are invertible matrices having the right triangular structure. In particular, this shows that $\Yc_2$ is invertible.\\
	\textbf{Step 2} (\textit{Formula for the triangular $\De_c$}).\\
	For reasons of space we drop the argument of $\De_c(.)$ and $\De_l(.)$. Motivated by the necessity part, the goal is to find a suitable triangular $\De_c$ structured as in (\ref{dec}) such that (\ref{cond_scheduling}) is true. We directly infer positive definiteness of the diagonal blocks in (\ref{cond_scheduling}) since $Q_2,\,Q_3\in\mathbf{P}_p$ and $\Qa_1\in\mathbf{P}_d$. Thus, an explicitly formula for $\De_c$ can be obtained by rendering in (\ref{cond_scheduling}) the off-diagonal blocks zero. Recall that $U_2$, $V_2$ are invertible, block-triangular matrices by construction which leads to the choice
	$$
	\De_c:=-U_2^{-T}
	\mat{cc}{Q_2\De_l\Qa_1+\De_l^T&0\\Q_3\De_l\Qa_1+\De_l^T&Q_3\De_l+\De_l^TQ_2}V_2^{-1}.
	$$
	By reversing the congruence transformation with $\Yc_2$ in the necessity part, (\ref{cond_scheduling}) implies (\ref{p31}) and thus $\Xc_2\in\mathbf{P}$.\\
	\textbf{Step 3} (\textit{Construction of controller matrices}).\\
	Let us now define $K_{ij},\,L_i,\,M_j$ by (\ref{s3}) and $N:=\NB=0$. Since $U_i,\,V_i$ are invertible for $i=1,2$, we can solve (\ref{s2}) for $A^c_{ij},\,B^c_i,\,C^c_j,\,D^c$; these controller matrices have indeed the desired structure of (\ref{s12}) as can be seen analogously to \cite{roesinger2019b} by exploiting the structure of $U_2,\,V_2$ and (\ref{rec}). Hence, (\ref{s17}) is true and, by applying congruence transformations with $\Yc_i^{-1}$ for $i=1,2$ along with the factorizations (\ref{fac}), we thus infer (\ref{s14}) for (\ref{s13}).~\hfill$\blacksquare$

\end{document}